\documentclass[11pt]{amsart}
\usepackage{amsmath,amssymb,latexsym}
\usepackage[T1]{fontenc}
\usepackage{enumerate}

\def\ale{\lesssim}
\def\FC{\mathrm{FC}}

\def\Ind#1#2{#1\setbox0=\hbox{$#1x$}\kern\wd0\hbox to 0pt{\hss$#1\mid$\hss}
\lower.9\ht0\hbox to 0pt{\hss$#1\smile$\hss}\kern\wd0}

\def\Notind#1#2{#1\setbox0=\hbox{$#1x$}\kern\wd0\hbox to 0pt{\mathchardef
\nn="3236\hss$#1\nn$\kern1.4\wd0\hss}\hbox to 0pt{\hss$#1\mid$\hss}\lower.9\ht0
\hbox to 0pt{\hss$#1\smile$\hss}\kern\wd0}

\theoremstyle{plain}
\newtheorem{theorem}{Theorem}[section]

\newtheorem{fact}[theorem]{Fact}

\theoremstyle{definition}
\newtheorem{defn}[theorem]{Definition}
\newtheorem{remark}[theorem]{Remark}

\def\pf{\par\noindent{\em Proof. }}

\title{A note on FC-nilpotency}
\author{Nadja Hempel}
\author{Daniel Palac\'in}
\thanks{The second author was partially supported by the project MTM2014--59178--P}
\address{Department of Mathematics, University of California Los Angeles, Los Angeles, CA 90095-1555, USA}
\email{nadja@math.ucla.edu}
\address{Einstein Institute of Mathematics, The Hebrew University of Jerusalem, Jerusalem, 9190401, Israel}
\email{daniel.palacin@mail.huji.ac.il}

\begin{document}

\begin{abstract}
The notion of bounded FC-nilpotent group is introduced and it is shown that any such group is nilpotent-by-finite, generalizing a result of Neumann on bounded FC-groups.
\end{abstract}

\maketitle

\section{Introduction}

A group in which every element has only finitely many conjugates is called a {\em finite conjugacy group}, FC-group for short. Of course, in particular all abelian groups and also all finite groups are FC-groups but there are many more non-trivial examples. The study of this class of groups was initiated by Baer \cite{Baer} and Neumann \cite{Neu}, and its general theory was strongly developed during the second half of the last century, see \cite{Tom}.

Numerous variations of the notion of FC-group have been considered to study structural properties of infinite groups with some finiteness condition. As a strengthening, Neumann considered FC-groups with a uniform bound on the size of the conjugacy classes, known as {\em bounded FC-groups}, and showed that these are finite-by-abelian \cite{Neu}. On the other hand, Haimo \cite{Haimo} and later Duguid and McLain
\cite{DM} analyzed FC-nilpotent and FC-solvable groups, which are natural generalizations of  nilpotent and solvable groups respectively to the FC context. For instance, a group is FC-solvable of length $n$ if it admits a finite chain of length $n$ of normal subgroups whose factors are FC-groups. Similarly, one can define the notion of FC-nilpotent, see Definition \ref{DefFCNil} for a precise definition. Furthermore, Hickin and Wenzel have shown in \cite{HicWen} that, like for normal nilpotent groups, the product of two normal FC-nilpotent subgroups is normal and again FC-nilpotent.

In this paper we aim to study a suitable version of bounded FC-nilpotency, see again Definition \ref{DefFCNil}, and to show that these groups are exactly the nilpotent-by-finite ones. This result generalizes the one of Neumann on bounded FC-groups and another of Duguid and McLain asserting that finitely generated FC-nilpotent groups are nilpotent-by-finite. Furthermore, bounded FC-nilpotent groups appear naturally in the study of groups in model theory such as in $\aleph_0$-categorical groups and groups definable in simple (or even wider families of) first-order theories. For instance, in these cases every definable FC-nilpotent group is indeed bounded. However, such groups are typically not finitely generated and therefore, the aforementioned result of Duguid and McLain cannot {\it a priori} be applied to deduce that these groups are nilpotent-by-finite. 

The result presented here generalizes some previous cases due to Wagner \cite[Proposition 4.4.10]{Wagner} for groups in simple theories, as well as in \cite{Nadja} for groups satisfying a uniform chain condition on centralizers up to bounded index. Our proof involves some machinery on FC-centralizers recently obtained by the first author in \cite{Nadja} using techniques from model theory. 
Finally, concerning the FC-solvable case note that the situation is more straightforward and an easy argument is given at the end of the paper.

\section{Bounded FC-nilpotent groups}

Given a group $G$ and a subset $X$ of $G$, we denote by $C_G(X)$ the elements of $G$ that commute with every element in $X$. Moreover, considering a normal subgroup $N$ of $G$, we denote by $C_G(g/N)$ the elements of $G$ in the preimage of $C_{G/N} (gN)$ under the usual projection.  

We recall the definition of an FC-centralizer due to Haimo \cite{Haimo} and related notions, which play an essential role along the paper. 

\begin{defn} A subgroup $H$ of $G$ is {\em contained up to finite index} in another subgroup $K$ if $H\cap K$ has finite index in $H$. We denote this by $H\ale K$. Then $H$ and $K$
are {\em commensurable}, denoted by $H\sim K$, if $H\ale K$ and $K\ale H$.
\end{defn}
Observe that $\ale$ is a transitive relation among subgroups of $G$, and that
$\sim$ is an equivalence relation.

\begin{defn}\label{DefFC}
Let $G$ be a group and let $K,H,N$ be subgroups of $G$ with $N$ normalized by $H$. The {\em FC-centralizer of $H$ modulo $N$ in $K$} is defined as
$$
\FC_K(H/N)=\{k\in N_K(N): H\sim C_H(k/N)\}.
$$
If $N$ is trivial it is omitted. %For a natural number $n$, the {\em $n$th iterated almost centralizer of $H$ in $G$} is defined inductively by
%$$
%\FC_G^1(H)=\FC_G(H) \ \mbox{ and } \ \FC_G^{n+1}(H)=\FC_G(H/\FC_G^n(H)).
%$$
%Moreover, we define the {\em almost center of $G$} as $\Z(G)=\FC_G(G)$ and the {\em $n$th iterated center of $G$} as $\Z_n(G)=\FC_G^n(G).$
\end{defn}
%For any subgroup $K$ of $G$ we put $\FC_K(H/N)=\FC_G(H/N)\cap K$. In the literature, the almost centralizer $\FC_K(H/N)$ is also denoted by by ${\rm FC}_K(H/N)$, and the almost center of a group $G$ as ${\rm FC}(G)$.
%However, we shall use the  terminology introduced above.

In other words, the group $\FC_K(H/N)$ consists of the elements $k$ in $N_K(N)$ such that $[H: C_H(k/N)]$ is finite, i.\ e.\  $k^H/N$ is finite.  Moreover, observe that $G$ is an FC-group if $G=\FC_G(G)$. As pointed out in the introduction a priori an FC-group may have arbitrarily large conjugacy classes. Those FC-groups in which there is a natural number bounding the size of any conjugacy class are called {\em bounded FC-groups}, and are precisely finite-by-abelian groups \cite[Theorem 5.1]{Neu}.

The following definition generalizes the notion of bounded FC-group to arbitrary FC-centralizers.

\begin{defn}
Let $G$ be a group and let $H,K,N$ be subgroups of $G$ with $K$ normalized by $H$. We say that $\FC_K(H/N)$ is {\em bounded} if there exists a natural number $n$ such that
$$
\FC_K(H/N)=\{g\in N_K(N): |H/C_H(g/N)|\le n\}.
$$
\end{defn}
%From a model-theoretic point of view, bounded almost centralizers are always definable subgroups and hence, the index of the centralizer of an element is preserved under taking elementary extensions. %In contrast, if an almost centralizer is not bounded, then model-theoretic compactness yields the existence of an elementary extension of the group in which
%This allows us to apply model-theoretic techniques to study bounded almost centralizers, up to finite index.
FC-centralizers for definable groups and bounded FC-centralizers have been studied by the first author in \cite{Nadja} who has shown that some behaviors of the ordinary centralizers can be generalized to FC-centralizers. This is exemplified in the following two lemmata which can be found as   \cite[Theorem 2.10]{Nadja} and \cite[Theorem 2.18]{Nadja}. For the former we give a proof which is a simple adaptation of the non-definable version to the bounded case.
\begin{fact}[Symmetry]\label{FactNadjaSym}
Let $H$ and $ K$  be two subgroups of a group $G$ and $N$ be a subgroup of $G$ that is normalized by $H$ and $K$. Assume further that the FC-centralizer of $\FC_H(K/N)$ is bounded. If $H\ale \FC_G(K/N)$, then $K\ale \FC_G(H/N)$.
\end{fact}

%The following proof is simple adaption of the proof given in \cite{Nadja} in the non-definable context.
\pf
Assume that $H$ and $\FC_H(K/N)$ be commensurable and let $d$ be the natural number such that  for any element $h$ in $\FC_H(K/N)$ we have that
$$
|H/C_H(h/N)| \leq d \ \ \ \ (\ast).
$$
Now suppose that $K$ and $\FC_K(H/N)$ are not commensurable. In this case, we can choose elements $k_0,\ldots,k_d$ in different cosets of $\FC_K(H/N)$ in $K$. Thus, the index $[H: C_H(k_i k_j^{-1}/N)]$ is infinite. Since no group can be covered by finitely many cosets of subgroups of infinite index by a well-known theorem of Neumann \cite{Neu2}, there are infinitely many elements $\{h_i\}_{i \in \mathbb N}$ in $H$ such that for all natural numbers $s\neq t$ and $ i\neq j \leq d$, we have that $[h_s h_t^{-1},k_i  k_j^{-1}] \not \in N$. Thus, the elements $k_0,\ldots, k_d$ witness that   $[K: C_H(h_s h_t^{-1}/N)] > d$. Hence, by $(\ast)$ none of the elements $h_s h_t^{-1}$ can belong to $\FC_H(K/N)$, and whence $H$ and $\FC_H(K/N)$ are not commensurable, contradicting our initial assumptions.
\qed
%Roughly speaking, this means that almost all elements of $H$ commute with almost all element of $K$ if and only if almost all elements of $K$ commutes with almost all elements of $K$. The saturation of our structure as well as the definability of $H$ and $K$, play an essential role. Furthermore,

\begin{fact}\label{FactNadjaCom}
Let $G$ be a group and let $H$ and $K$ be two subgroups of $G$ such that $H$ is normalized by $K$. If $H=\FC_H(K)$, $K=\FC_K(H)$, and $\FC_K(H)$ is bounded, then the commutator subgroup $[H,K]$ is finite.
\end{fact}
Notice that the latter generalizes the aforementioned result due to Neumann on bounded FC-groups.

\begin{defn}\label{DefFCNil}
A group $H$ is {\em FC-nilpotent} of class $n$ if there exists a finite chain of normal subgroups
$$
\{1\}\le H_1 \le H_2 \le \dots \le H_n=H
$$
such that $H_{i+1}$ is contained in $\FC_H(H/H_i)$. We say that $H$ is {\em bounded FC-nilpotent} if additionally each $\FC_{H_{i+1}}(H/H_i)$ is bounded.
\end{defn}

FC-nilpotent groups were introduced by Haimo \cite{Haimo}, who studied their basic properties. Duguid and McLain \cite[Corollary 1]{DM} proved that finitely generated FC-nilpotent groups of class $n$ are (nilpotent of  class $n$)-by-finite.

\begin{remark} A nilpotent-by-finite group is bounded FC-nilpotent.\end{remark}
\pf Let $H$ be a nilpotent-by-finite group and choose $N$ a normal nilpotent finite index subgroup of $H$. Let $n$ be the nilpotency class of $N$ and $k$ be the index of $N$ in $H$. It is clear that $Z_{i+1}(N)$ is contained in $\FC_H(H/Z_i(N))$ and that
$$
\big[H:C_H(x/Z_i(N))\big]\le k
$$
for any $x\in Z_{i+1}(N)$. Then
$$
\{1\} \le Z(N) \le \dots \le Z_{n}(N) \le H
$$
witnesses that $H$ is a bounded FC-nilpotent group.
 \qed

The following result yields the equivalence between bounded FC-nilpotent groups and nilpotent-by-finite ones.

\begin{theorem}
Any bounded FC-nilpotent group of class $n$ is (nilpotent of class $2n$)-by-finite.
\end{theorem}
\pf
Let $N$ be a bounded FC-nilpotent group of class $n$, and let
$$
\{1\}=N_0\le N_1 \le \dots \le N_n=N
$$
be a series witnessing this. %Consider the first-order structure of $N$ as a pure group together with predicates for the subgroups $N_i$ for $i\le n$. As each almost centralizer $\FC_{N_{i+1}}(N/N_i)$ is bounded, it is definable in this structure and so, in any elementary extension the corresponding subgroups $\bar N_0,\ldots,\bar N_n$ satisfy that each $\bar N_{i+1}$ equals to $\FC_{\bar N_{i+1}}(\bar N/\bar N_i)$ and the latter is bounded. Thus,  the subgroups $\bar N_0,\ldots,\bar N_n$ witness that $\bar N$ is bounded FC-nilpotent of class $n$ as well. Hence, without loss we may assume that the structure $(N,\cdot,N_0,\ldots,N_n)$ is $\aleph_0$-saturated and work inside it. From now on, by definable we mean a definable object in the expanded structure.
Now, we find recursively on $i\le n$ a subgroup $F_i$ of finite index in $N$ and a finite chain
$$
\{1\}=H_0^i \le H_1^i \le \dots \le H_{2i}^i
$$ of normal subgroups of $F_i$ such that $H_{j+1}^i/H_j^i$ is central in $F_i/H_j^i$ for $j<2i$ and $F_i \cap N_i$ is contained in $H_{2i}^i$.
Once this process is concluded, the subgroup $F_n$ is nilpotent of class at most $2n$ and has finite index in $N$.

We start the construction by setting $H_0^0=\{1\}$ and $F_0=N$. To continue, let $i>0$ and assume that we have already defined $H_0^{i-1},\ldots,H_{2(i-1)}^{i-1}$ and $F_{i-1}$.
%As $N_i$ is contained in $\FC_N(N/N_{i-1})$ and $\FC_{N_i}(N/N_{i-1})$ is bounded, then $N_i\cap F_{i-1}$ is also contained in $\FC_{F_{i-1}}(F_{i-1}/F_{i-1}\cap N_{i-1})$ and the latter is bounded. Hence, there is no harm to assume that $N=F_{i-1}$ and $N_i=F_{i-1} \cap N_i$.
Now, we work in the quotient $\bar N = F_{i-1}/ H_{2(i-1)}^{i-1}$ and consider the subgroup $\bar N_i = N_i \cap F_{i-1}/ H_{2(i-1)}^{i-1}$. As $N_{i-1}\cap F_{i-1}$ is contained in $H_{2(i-1)}^{i-1}$ and $N_i=\FC_{N_i}(N/N_{i-1})$ is bounded, we have that $\bar N_i =\FC_{\bar N_i }(\bar N)$ and that it is bounded.  Thus Fact \ref{FactNadjaSym} yields that $\bar N\sim \FC_{\bar N}(\bar N_i)$. % Despite that $\FC_{\bar N}(\bar N_i)$ is not necessarily definable, it is given as the union of definable sets
%$$
%\FC_{\bar N} \big( \bar N_i \big)_k = \big\{x\in \bar N : \big|\mathcal  N_i/C_{\bar N_i}(x) \big| \le k \big\}.
%$$
%Thus, there are finitely many elements $x_1,\ldots,x_m$ in $\bar N$ such that
%$$
% \bar N=\bigcup_{i=1}^m x_i\cdot \Big( \bigcup_{k\in\mathbb N} \FC_{ \bar N }\big( \bar N_i \big)_k \Big) = \bigcup_{i=1}^m \bigcup_{k\in\mathbb N} x_i\cdot\FC_{\bar N} \big( \bar N_i \big)_k
%$$
%and thus since $\bar N$ is definable, model-theoretic compactness yields that
%$$
%\bar N = \bigcup_{i=1}^m x_i\cdot \FC_{\bar N} \big(\bar N_i \big)_k
%$$
%for some $k$. Hence $\FC_{\bar N} (\bar N_i)$ equals $\FC_{\bar N}( \bar N_i)_k$,  so it is bounded and whence definable.

Let $\bar N^\ast$ be equal to $\FC_{\bar N} (\bar N_i)$, a subgroup of finite index in $\bar N$, and let $\bar N_i^\ast$ be equal to $\bar N_i\cap \bar N^\ast$. Now, as $\bar N_i^\ast$ has also finite index in $\bar N_i$, we clearly have that $\bar N^\ast = \FC_{\bar N^\ast} (\bar N_i^\ast)$
and moreover
$$
\bar N_i ^\ast= \FC_{\bar N_i }(\bar N) \cap \bar N^\ast = \FC_{\bar N_i^\ast }(\bar N^\ast).
$$
As the latter remains bounded,  Fact \ref{FactNadjaCom} yields that the group $X$ defined as $[\bar N^\ast, \bar N_i^\ast]$ is finite. Note that all  considered groups are normal in the ambient group. Hence $X$ is contained in $\bar N_i^\ast= \FC_{\bar N_i^\ast}(\bar N^\ast)$ and whence $C_{\bar N^\ast}(X)$ has finite index in $\bar N^\ast$.

%Let $G$ be $\FC_N( N_i/N_{i-1})$, a bounded subgroup of finite index in $N$. Thus $N_i/N_{i-1}$ equals the bounded group $\FC_{N_i/N_{i-1}}(G/N_{i-1})$ and so the commutator subgroup $[N_i,G]$ is finite modulo $N_{i-1}$ by Fact \ref{FactNadjaCom} as $G/N_{i-1}$ is trivially contained in $\FC_{G/N_{i-1}}(N_i/N_{i-1})$ by definition.

Finally, set $F_i$ to be the preimage of $C_{\bar N^\ast}(X)$ in $N$, which as well has finite index in $N$, the group  $H_{2i-1}^i$  to be to preimage of $X\cap C_{\bar N^\ast}(X) $ in $N$, and $H_{2i}^i$ to be the preimage of $N_i^\ast \cap C_{\bar N^\ast}(X)$ in $N$.  Moreover,  for $j < 2i-1$, we let $H_{j}^i$ be  $H_{j}^{i-1} \cap F_i$. Note that by construction, the subgroups $H^i_{j}$ for $j\leq 2i$ are contained in $F_i$. Then the sequence
$$
\{1\}=H_0^{i} \le H_1^{i} \le \dots \le H_{2i}^{i}
$$
together with $F_i$ is as desired. \qed

%To conclude, set $H_j^i=F_i\cap H_j^{i-1}$ for $j\le 2(i-1)$,  $H_{2i-1}^i=C_{F_i}(F_i/F_i\cap N_{i-1})$ and $H_{2i}^i=C_{F_i}^2(F_i/F_i\cap N_{i-1})$. It is routine to check that this subgroups are as desired. Indeed, for $j\le 2(i-1)$ by construction we have that  $H_{j+1}^i/H_j^i$ is central in $F_i/H_j^i$ and that $F_i\cap N_{i-1}$ is clearly contained in $H_{2i-2}^i$. Thus, the latter yields that $H_{2i-1}^i$ and $H_{2i}^i$ are contained in $C_{F_i}(F_i/H_{2i-2}^i)$ and in $C_{F_i}^2(F_i/H_{2i-2}^i)$ respectively. That is, the groups $H_{2i-1}^i/H_{2i-2}^i$ and $H_{2i}^i/H_{2i-1}^i$ are central in $F_i/H_{2i-2}^i$ and in $F_i/H_{2i-1}^i$ respectively. Finally, as
%$$
%[F_i\cap N_i,F_i] \le [N_i,G]
%$$
%we get that $[F_i\cap N_i ,_2 F_i]$ is trivial by definition of $F_i$ and so $F_i\cap N_i$ is contained in $H_{2i}^i$. This finishes the proof.
%\qed

\begin{remark}
Observe that our proof yields that each quotient $H_{2i+1}^i$ modulo $H_{2i}^i$ is finite. Thus, the group $N$ has a finite index nilpotent subgroup of class $2n$, which admits a series of length $n$ where each factor is (finite central)-by-central.
\end{remark}

\section{Bounded FC-solvability}
We say that a group is {\em bounded FC-solvable} of length $n$ if it is an FC-solvable group of length $n$ in which each factor witnessing the FC-solvability is a bounded FC-group. Thus, by Neumann's Theorem each of these factors is finite-by-abelian. Moreover, for any finite-by-abelian group $H$, the characteristic subgroup $C_H(H')$ is nilpotent of class two and has finite index. Thus, such a group is nilpotent-by-finite. Using this, we can easily show by induction on the FC-solvability length that a bounded FC-solvable group is solvable-by-finite:

More precisely, let 
$$
\{ 1 \} = G_0 \unlhd G_1 \unlhd \ldots \unlhd G_{n-1} \ldots \unlhd G_n=G
$$
be a series that witnesses that $G$ is a bounded FC-solvable group. Then, the group $G_{n-1}$ is FC-solvable of smaller length and hence solvable-by-finite. Thus $G$ is a solvable-by-finite-by-abelian group. So, as any   finite-by-abelian group is a nilpotent-by-finite group, we can conclude that $G$ is solvable-by-finite.

% is nilpotent of class two ab a. Then $G^\ast=C_G(G_1')$ is a finite index subgroup of $G$ such that  and then the group $G^\ast/Z()$ is FC-solvable of smaller length and hence solvable-by-finite. Whence $G$ itself is solvable-by-finite since $C_{G_1}(G_1')$ is nilpotent of class two.

\end{document}